\documentclass[10pt,sumlimits,namelimits]{article}
\newcommand{\nc}{\newcommand}
\newcommand{\rnc}{\renewcommand}
\usepackage{amsthm}
\usepackage{amsmath}
\usepackage{amsfonts}
\usepackage{amssymb}
\usepackage{color}
\usepackage[colorlinks=true,urlcolor=blue]{hyperref}
\usepackage[applemac]{inputenc}

\nc{\N}{\mathbb{N}}
\nc{\Z}{\mathbb{Z}}
\nc{\D}{\mathbb{D}}
\nc{\Q}{\mathbb{Q}}
\nc{\R}{\mathbb{R}}
\nc{\C}{\mathbb{C}}
\nc{\vD}{{\cal D}}
\nc{\vF}{{\cal F}}
\nc{\vS}{{\cal S}}
\nc{\vphi}{\varphi}
\nc{\eps}{\varepsilon}
\nc{\dsp}{\displaystyle}
\nc{\ovl}{\overline}
\nc{\udl}{\underline}
\nc{\vlim}{\lim\limits}
\nc{\vlimsup}{\limsup\limits}
\nc{\vliminf}{\liminf\limits}
\nc{\vsup}{\sup\limits}
\nc{\vinf}{\inf\limits}
\nc{\vint}{\int\limits}
\nc{\inj}{\hookrightarrow}
\nc{\tends}{\longrightarrow}
\nc{\weak}{\rightharpoonup}
\nc{\w}{{\textsl w}}
\nc{\loc}{{\rm loc}}
\rnc{\b}{{\rm b}}
\rnc{\le}{\leqslant}
\rnc{\ge}{\geqslant}
\rnc{\Re}{{\rm Re}}
\rnc{\Im}{{\rm Im}}

\numberwithin{equation}{section}
\rnc{\theequation}{\thesection.\arabic{equation}}

\newtheorem{thm}{Theorem}[section]

\theoremstyle{definition}
\newtheorem{rmk}[thm]{Remark}
\newtheorem{defi}[thm]{Definition}
\newtheorem{exa}[thm]{Example}

\newenvironment{proof*}{\noindent{\bf Proof.}}{\qed}
\newenvironment{vproof}[1]{\noindent{\bf Proof #1}}{\qed}

\title{\Huge \sc Maximum Decay Rate for Finite--Energy Solutions of Nonlinear Schrödinger Equations}
\author{\sc Pascal Bégout}
\date{}

\begin{document}

\maketitle

\begin{center}
Laboratoire Jacques-Louis Lions \\
Université Pierre et Marie Curie \\
Boîte Courrier 187 \\
4, place Jussieu 75252 Paris Cedex 05, FRANCE
\bigskip \\
{\footnotesize e-mail\:: }\htmladdnormallink{{\footnotesize\udl{\tt{begout@ann.jussieu.fr}}}}
{mailto:begout@ann.jussieu.fr}
\end{center}

\begin{abstract}
We give explicit time lower bounds in the Lebesgue spaces for all nontrivial solutions of nonlinear Schrödinger equations bounded in the energy space. The result applies for these equations set in any domain of $\R^N,$ including the whole space. This also holds for a large class of nonlinearities, thereby extending the results obtained by Hayashi and Ozawa in~\cite{MR91d:35035} and by the author in~\cite{beg3}.
\end{abstract}

\baselineskip .7cm

\section{Introduction and notations}
\label{introduction}

{\let\thefootnote\relax\footnotetext{2000 Mathematics Subject Classification: 35Q55,(35B40)}}

We consider global solutions of the following nonlinear Schrödinger equation,
\begin{gather}
 \left\{
  \begin{split}
   \label{nls}
    i\frac{\partial u}{\partial t}+\Delta u+f(u)= 0, & \; (t,x)\in\R\times\Omega, \\
                                 	u_{|\partial\Omega}= 0, & \mbox{ in }\R\times\partial\Omega, \\
	                                                   u(0)= u_0, & \mbox{ in }\Omega, \\
  \end{split}
 \right.
\end{gather}
where $\Omega\subseteq\R^N$ is an open subset with boundary $\partial\Omega,$ the nonlinearity $f$ satisfies some conditions to be specified later and $u_0$ is a given initial data. We show that if $u$ is a solution of (\ref{nls}) with initial data $u_0\not\equiv0,$ then for any $0<\eps<\|u_0\|_{L^2(\Omega)},$ $r\in[2,\infty]$ and $c\in\Omega,$ we have
\begin{gather}
\label{intro0}
\liminf_{t\to\pm\infty}|t|^{N\left(\frac{1}{2}-\frac{1}{r}\right)}\|u(t)\|_{L^r(\Omega\cap\{|x-c|<M_0|t|\})}
\ge\eps|B(0,M_0)|^{-\left(\frac{1}{2}-\frac{1}{r}\right)},
\end{gather}
where $M_0>0$ is an explicit constant depending only on $\|u\|_{L^\infty(\R;H^1_0(\Omega))}$ and $\eps.$ \\
In \cite{297.35062}, Strauss treated the case of the free operator. He showed that for any $u_0\in L^2(\R^N)\setminus\{0\}$ and any $r\in[2,\infty],$
\begin{gather}
\label{intro1}
\liminf_{t\to\pm\infty}|t|^{N\left(\frac{1}{2}-\frac{1}{r}\right)}\|u(t)\|_{L^r(\Omega_t)}>0,
\end{gather}
where $u(t)=e^{it\Delta}u_0$ and $\Omega_t=\{|x|<k|t|\},$ for some $k>0$ (see Strauss~\cite{297.35062}, Lemma~p.69). If fact, as seen in Ozawa~\cite{MR91m:35183}, the same result still holds with $\Omega_t=\{k'|t|<|x|<k|t|\},$ for some $k>k'>0.$ Later, (\ref{intro1}) was  established with $\Omega_t\equiv\R^N$ and still for the free operator, but for any $u_0\in\vS'(\R^N)\setminus\{0\}$ and any $r\in[1,\infty]$ (see Kato~\cite{MR95i:35276}, Decay Lemma~p.228). Ozawa~\cite{MR91m:35183} showed that, for certain potentials $V:\R^N\tends\R,$ any nontrivial asymptotically free solution $u(t)=e^{it(\Delta+V)}u_0$ satisfies (\ref{intro1}), with $\Omega_t=\{k'|t|<|x|<k|t|\}$ for some $k>k'>0.$ By {\it asymptotically free}, we mean that $\|u(t)-e^{it\Delta}u_+\|_{L^2(\R^N)}\xrightarrow{t\to\infty}0,$ for some $u_+\in L^2(\R^N).$ However, the constants $k$ and $k'$ are not explicit. \\
The nonlinear case was first treated by Hayashi and Ozawa~\cite{MR91d:35035}, still with $\Omega=\R^N.$  They showed that any nontrivial solution $u$ of (\ref{nls}) with initial data $u_0\in H^1(\R^N)$ satisfies (\ref{intro1}), where $\Omega_t=\{|x|<k|t|\}$ for some $k>0.$ The nonlinearity $f$ had to be a single power interaction $(f(u)=-|u|^\alpha u),$ nonlocal interaction $(f(u)=-(|x|^{-\mu}*|u|^2)u)$ or external potential $(f(u)=Vu),$ and in all those cases, repulsive. Furthermore, the initial data had to satisfy the additional assumption $\int_{\R^N}|x|^2|u_0(x)|^2dx<\infty.$ These two last assumptions were needed in order to use the pseudo-conformal transformation law and obtain some {\it a priori} estimates on the time decay of the solution. \\
Finally, in \cite{beg3} the author considered the special case where the nonlinearity is a single power interaction, with $\Omega=\Omega_t=\R^N.$ He established (\ref{intro1}) for all nontrivial solutions $u$ of (\ref{nls}) with initial value in $H^1(\R^N).$ In particular, estimate (\ref{intro1}) still holds for the attractive nonlinearity $f(u)=|u|^\alpha u.$ \\
Note that all the above results are obtained by contradiction, which explains the fact that the lower bound in (\ref{intro1}) is not explicit. \\
In this paper, we establish estimate (\ref{intro0}) with a very simple (and direct) proof. It only makes use the linear part of equation (\ref{nls}). This permits us to extend the results of Hayashi and Ozawa in \cite{MR91d:35035} to any domain $\Omega\subseteq\R^N$ and to a large class of nonlinearities. The proof requires that for almost every $u\in\C,$ $\Im (f(u)\ovl u)=0$ and that the solution is bounded in the energy space, namely in $H^1_0(\Omega).$ Note that the first hypothesis is quite reasonable. It is needed to have conservation of charge, which is essential to obtain (\ref{intro0}). The proof allows us to consider attractive nonlinearities. In particular, we prove (\ref{intro0}) for the free operator with $\Omega\not=\R^N.$ Finally, it seems that the possibility of choosing the lower bound in the right hand side of (\ref{intro0}) for the solutions of nonlinear Schrödinger equations (\ref{nls}) is new. \\
This paper is organized as follows. In Section~\ref{slb}, we give and establish the main result, namely the estimate (\ref{intro0}). In Section~\ref{example}, we give some examples for which the result of Section~\ref{slb} applies. \\
Throughout this paper, we use the following notation. We denote by $\ovl z$ the conjugate of the complex number $z,$ by $\Re(z)$ its real part and by $\Im(z)$ its imaginary part. We denote by $\Omega\subseteq\R^N$ any nonempty open subset. For $1\le p\le\infty,$ $p'$ denotes the conjugate of $p$ defined by $\frac{1}{p}+\frac{1}{p'}=1;$ $L^p(\Omega)=L^p(\Omega;\C)$ is the usual Lebesgue space and we write $\|f\|_{L^p(\Omega)}=\infty$ if $f\in L^1_\loc(\Omega)$ and $f\not\in L^p(\Omega).$ $H^1_0(\Omega)=H^1_0(\Omega;\C)$ is the usual Sobolev space and $H^{-1}(\Omega)$ is its topological dual. The Laplacian in $\Omega$ is written $\Delta=\sum\limits_{j=1}^N\frac{\partial^2}{\partial x^2_j}.$ For a Banach space $(E,\| \: . \: \|_E),$ we denote by $E^*$ its topological dual and by $\langle\: . \; , \: . \:\rangle_{E^*,E}\in\R$ the $E^*-E$ duality product. In particular, for any $T\in L^{p'}(\Omega)$ and $\vphi\in L^p(\Omega)$ with $1\le p<\infty,$ $\langle T,\vphi\rangle_{L^{p'}(\Omega),L^p(\Omega)}=\Re\vint_\Omega T(x)\ovl{\vphi(x)}dx.$ Finally, for $c\in E$ and $R\in(0,\infty),$ we denote by $B_E(c,R)=\{x\in E;\; \|x-c\|_E<R\}$ the open ball of $E$ of center $c$ and radius $R.$

\section{Sharp lower bound}
\label{slb}

Let $\Omega\subseteq\R^N$ be an open subset and $f$ satisfying the following assumptions.
\begin{gather}
 \label{hypf1}
  f\in C(H^1_0(\Omega);H^{-1}(\Omega)), \\
 \label{hypf2}
  \forall a\in W^{1,\infty}(\Omega;\R),\; \forall u\in H^1_0(\Omega;\C),\;
  \langle f(u),iau\rangle_{H^{-1}(\Omega),H^1_0(\Omega)}=0.
\end{gather}
When the solution is not smooth enough, the nonlinearity $f$ has to satisfy the additional assumption
\begin{gather}
 \label{hypf3}
  f:H^1_0(\Omega)\tends H^{-1}(\Omega) \mbox{ is bounded on bounded sets.}
\end{gather}

\begin{defi}
\label{defisol}
Let $\Omega\subseteq\R^N$ be an open subset, let $f$ satisfying (\ref{hypf1}), let $u_0\in H^1_0(\Omega)$ and let $u\in L^\infty((0,\infty);H^1_0(\Omega)).$ Assume further that for any $T>0,$ $f(u)\in L^\infty((0,T);H^{-1}(\Omega)).$ Then we say that $u$ is a solution of (\ref{nls}) if $u(0)=u_0$ and if for any $\vphi\in H^1_0(\Omega;\C)$ and $\psi\in\vD((0,\infty);\R),$
\begin{gather}
\label{edpsol}
\vint_0^{+\infty}\left\{-\langle iu(t),\vphi\rangle_{L^2(\Omega),L^2(\Omega)}\psi'(t)
+\langle\Delta u(t)+f(u(t)),\vphi\rangle_{H^{-1}(\Omega),H^1_0(\Omega)}\psi(t)\right\}dt=0.
\end{gather}
In the same way, we define a solution $u\in L^\infty((-\infty,0);H^1_0(\Omega))$ and $u\in L^\infty(\R;H^1_0(\Omega)).$
\end{defi}

\begin{rmk}
\label{rmkmaxdecaylinear8}
Note that Definition \ref{defisol} makes sense since we then have $u\in W^{1,\infty}((0,T);H^{-1}(\Omega)),$ for any $T>0$ and so the first equation in (\ref{nls}) makes sense in $H^{-1}(\Omega)$ for almost every $t>0.$ Moreover, it follows from the inequality $\|v\|_{L^2(\Omega)}^2\le\|v\|_{H^{-1}(\Omega)}\|v\|_{H^1_0(\Omega)}$ that $u\in C([0,\infty);L^2(\Omega))$ and so $u(0)=u_0$ takes sense in $L^2(\Omega).$
\end{rmk}
The main result of this section is the following.

\begin{thm}
\label{thmmaxdecaylinear}
Let $\Omega\subseteq\R^N$ be an open subset, let $f$ satisfying $(\ref{hypf1})-(\ref{hypf2}),$ let $u_0\in H^1_0(\Omega)\setminus\{0\}$ and let $u\in L^\infty((0,\infty);H^1_0(\Omega)).$ Assume further that for any $T>0,$ $f(u)\in L^\infty((0,T);H^{-1}(\Omega))$ and that $u$ is a solution of $(\ref{nls})$ with $u(0)=u_0$ $($Definition~$\ref{defisol}).$ Then the following holds. Let $0<\eps<\|u_0\|_{L^2(\Omega)}$ and let
\begin{gather}
\label{thmmaxdecaylinear0}
M_0=\frac{2\|u_0\|_{L^2(\Omega)}\|\nabla u\|_{L^\infty((0,\infty);L^2(\Omega))}}{\|u_0\|_{L^2(\Omega)}^2-\eps^2}.
\end{gather}
Then we have for any $ r\in[2,\infty]$ and any $c\in\Omega,$
\begin{gather}
 \label{thmmaxdecaylinear2}
  \liminf_{t\to\infty}|t|^{N\left(\frac{1}{2}-\frac{1}{r}\right)}\|u(t)\|_{L^r(\Omega\cap\{|x-c|<M_0|t|\})}
  \ge\dfrac{\eps}{|B(0,M_0)|^{\frac{1}{2}-\frac{1}{r}}}.
\end{gather}
A similar statement holds for $t<0,$ with the obvious modification in $(\ref{thmmaxdecaylinear0}).$
\end{thm}

\begin{rmk}
\label{rmkmaxdecaylinear6}
Let $u\in L^\infty((0,\infty);H^1_0(\Omega))$ and $f\in C(H^1_0(\Omega);H^{-1}(\Omega)).$ If, in addition, $f:H^1_0(\Omega)\tends H^{-1}(\Omega)$ is bounded on bounded sets or if $u\in C([0,\infty);H^1_0(\Omega))$ then we obviously have that  for any $T>0,$ $f(u)\in L^\infty((0,T);H^{-1}(\Omega)).$
\end{rmk}

\begin{rmk}
\label{rmkmaxdecaylinear7}
Let $u$ be a solution of (\ref{nls}) (Definition \ref{defisol}).  It follows from Remark~\ref{rmkmaxdecaylinear8} that  for any $t\ge0,$ $u(t)\in L^2(\Omega).$ Assume further that $f$ satisfies (\ref{hypf2}). Then it follows immediately from $(\ref{hypf2})$ that conservation of charge holds, that is for any $t\ge0,$ $\|u(t)\|_{L^2(\Omega)}=\|u_0\|_{L^2(\Omega)}.$
\end{rmk}

\begin{rmk}
\label{rmkmaxdecaylinear5}
Notice that uniqueness of the solution is not required.
\end{rmk}

\begin{rmk}
\label{rmkmaxdecaylinear1}
Estimate (\ref{thmmaxdecaylinear2}) is trivial if $\Omega$ is bounded. Indeed, we have in this case that $\Omega\cap B(c,M|t|)=\Omega,$ for any $c\in\Omega,$ $M>0$ and $t\in\R$ large enough. Therefore, from conservation of charge and Hölder's inequality, we have
$$
\|u(t)\|_{L^r(\Omega)}\ge|\Omega|^{-\left(\frac{1}{2}-\frac{1}{r}\right)}\|u(t)\|_{L^2(\Omega)}
=|\Omega|^{-\left(\frac{1}{2}-\frac{1}{r}\right)}\|u_0\|_{L^2(\Omega)},
$$
for any $r\in[2,\infty]$ and any $t.$
\end{rmk}

\begin{rmk}
\label{rmkmaxdecaylinear9}
Note that in the particular case where $r=2,$ estimate (\ref{thmmaxdecaylinear2}) becomes
$$
\forall c\in\Omega,\; \liminf_{t\to\infty}\|u(t)\|_{L^2(\Omega\cap\{|x-c|<M_0|t|\})}\ge\eps.
$$
\end{rmk}

\begin{vproof}{of Theorem \ref{thmmaxdecaylinear}.}
We follow a method of Cazenave~\cite{caz1,MR2002047} (see Step 1 of the proof of Theorem~7.5.1 in~\cite{caz1} or Step 1 of the proof of Theorem~7.7.1 in~\cite{MR2002047}).
We write $I=(0,\infty)$ if $u$ is positively global, and $I=(-\infty,0)$ if $u$ is negatively global.
Let $\eps\in(0,\|u_0\|_{L^2(\Omega)}),$ let $M_0$ be given by (\ref{thmmaxdecaylinear0}) and let $c\in\Omega.$
For $H>0,$ define $T_H\in W^{1,\infty}(\Omega;\R)$ by
\begin{gather*}
\forall x\in\Omega,\; T_H(x)=
 \begin{cases}
  \: 1-\dfrac{|x-c|}{H}, & \mbox{ if } |x-c|<H, \medskip \\
  \:                            0, & \mbox{ if } |x-c|\ge H.
 \end{cases}
\mbox{ Then }\|T_H\|_{L^\infty}=1, \: \|\nabla T_H\|_{L^\infty}=\frac{1}{H} \: \mbox{and it follows }
\end{gather*}
that $iT_Hu\in L^\infty(I,H^1_0(\Omega)).$ By Remark~\ref{rmkmaxdecaylinear8}, the first equation in (\ref{nls}) makes sense in $H^{-1}(\Omega)$ for almost every $t\in I$ and so we can take the $H^{-1}-H^1_0$ duality product with $iT_Hu.$ Thus with (\ref{hypf2}), this yields
\begin{gather}
\label{proofthmmaxdecaylinear1}
\langle iu_t+\Delta u+f(u),iT_Hu\rangle_{H^{-1}(\Omega),H^1_0(\Omega)}=0, \medskip \\
\label{proofthmmaxdecaylinear2}
\langle iu_t,iT_Hu\rangle_{H^{-1}(\Omega),H^1_0(\Omega)}=\frac{1}{2}\frac{d}{dt}\vint_\Omega T_H(x)|u(\: . \:,x)|^2dx, \medskip \\
\label{proofthmmaxdecaylinear3}
\langle f(u),iT_Hu\rangle_{H^{-1}(\Omega),H^1_0(\Omega)}=0, \medskip \\
\label{proofthmmaxdecaylinear4}
\langle\Delta u,iT_Hu\rangle_{H^{-1}(\Omega),H^1_0(\Omega)}
=-\langle\nabla u,iu\nabla T_H\rangle_{L^2(\Omega),L^2(\Omega)}
=-\Im\vint_\Omega\ovl{u(\: . \:,x)}\nabla u(\: . \:,x).\nabla T_H(x)dx,
\end{gather}
almost everywhere on $I.$ It follows from (\ref{proofthmmaxdecaylinear4}), Hölder's inequality and conservation of charge (Remark~\ref{rmkmaxdecaylinear7}), that
\begin{gather}
\label{proofthmmaxdecaylinear5}
\left|\langle\Delta u,iT_Hu\rangle_{H^{-1}(\Omega),H^1_0(\Omega)}\right|
\le\dfrac{1}{H}\|u_0\|_{L^2(\Omega)}\|\nabla u\|_{L^\infty(I,L^2(\Omega))},
\end{gather}
almost everywhere on $I.$ Let $t\in I.$ Putting together (\ref{proofthmmaxdecaylinear1})--(\ref{proofthmmaxdecaylinear4}), integrating over $(0,t)$ if $t>0,$ or over $(t,0)$ if $t<0$ and using (\ref{proofthmmaxdecaylinear5}), we see that
$$
\vint_\Omega T_H(x)|u(t,x)|^2dx\ge\vint_\Omega T_H(x)|u_0(x)|^2dx-\dfrac{2|t|}{H}\|u_0\|_{L^2(\Omega)}\|\nabla u\|_{L^\infty(I,L^2(\Omega))}.
$$
We choose $H=M_0|t|$ in the above estimate. This gives
$$
\forall t\in I,\; \|u(t)\|^2_{L^2(\Omega\cap\{|x-c|<M_0|t|\})}\ge\vint_\Omega T_{M_0|t|}(x)|u_0(x)|^2dx-\dfrac{2}{M_0}\|u_0\|_{L^2(\Omega)}\|\nabla u\|_{L^\infty(I,L^2(\Omega))}.
$$
By the dominated convergence Theorem, we have $\dsp\vlim_{t\to\pm\infty}\vint_\Omega T_{M_0|t|}(x)|u_0(x)|^2dx=\|u_0\|^2_{L^2(\Omega)}$ which yields with the above estimate
\begin{gather}
\label{proofthmmaxdecaylinear6}
\liminf_{\underset{t\in I}{|t|\to\infty}}\|u(t)\|^2_{L^2(\Omega\cap\{|x-c|<M_0|t|\})}\ge\|u_0\|^2_{L^2(\Omega)}-\dfrac{2}{M_0}\|u_0\|_{L^2(\Omega)}\|\nabla u\|_{L^\infty(I,L^2(\Omega))}\stackrel{{\rm def}}{=}\eps^2.
\end{gather}
Let $r\in(2,\infty].$ It follows from Hölder's inequality that for any $t\in I,$
\begin{align*}
\|u(t)\|_{L^2(\Omega\cap\{|x-c|<M_0|t|\})}\le & |B(c,M_0|t|)|^{\frac{1}{2}-\frac{1}{r}}\|u(t)\|_{L^r(\Omega\cap\{|x-c|<M_0|t|\})} \\
= & |B(0,M_0)|^{\frac{1}{2}-\frac{1}{r}}|t|^{N\left(\frac{1}{2}-\frac{1}{r}\right)}\|u(t)\|_{L^r(\Omega\cap\{|x-c|<M_0|t|\})}.
\end{align*}
From the above estimate and from (\ref{proofthmmaxdecaylinear6}), it follows that
\begin{gather*}
\liminf_{\underset{t\in I}{|t|\to\infty}}|t|^{N\left(\frac{1}{2}-\frac{1}{r}\right)}\|u(t)\|_{L^r(\Omega\cap\{|x-c|<M_0|t|\})}
\ge\eps|B(0,M_0)|^{-\left(\frac{1}{2}-\frac{1}{r}\right)},
\end{gather*}
which is (\ref{thmmaxdecaylinear2}). Hence the result.
\end{vproof}

\begin{rmk}
\label{rmkmaxdecaylinear3}
In some cases, it may happen that $f$ does not satisfy (\ref{hypf1})--(\ref{hypf3}), but some slightly different assumptions. Let $\vD(\Omega)\subset E\inj H^1_0(\Omega)$ be a Banach space such that for any $a\in W^{1,\infty}(\Omega;\R)$ and $u\in E,$ $au\in E.$ Assume that $f$ satisfies the following assumptions.  
\begin{gather}
 \label{hypf1'}
  f\in C(E;E^*), \\
 \label{hypf2'}
  \forall a\in W^{1,\infty}(\Omega;\R),\; \forall u\in E,\;
  \langle f(u),iau\rangle_{E^*,E}=0, \\
  \label{hypf3'}
  \forall A>0,\; \exists R>0\; \mbox{ such that }\; f(B_E(0,A))\subset B_{E^*}(0,R).
\end{gather}
Let $u_0\in E\setminus\{0\}$ and assume that there exists a solution $u\in L^\infty((0,\infty);E)$ of (\ref{nls}) with initial data $u_0.$ We claim that conclusion of Theorem~\ref{thmmaxdecaylinear} holds. Indeed,
from (\ref{nls}), (\ref{hypf1'}) and (\ref{hypf3'}), we have that $u_t\in L^\infty((0,\infty);E^*).$ Therefore, the first equation in $(\ref{nls})$ makes sense in $E^*$ for almost every $t>0,$ and inequality
$\|u\|_{L^2(\Omega)}^2\le\|u\|_{E^*}\|u\|_E$ implies that $u\in C([0,\infty);L^2(\Omega)).$ Thus $u(0)=u_0$ takes sense in $L^2(\Omega)$ and we can take  the $E^*-E$ duality product of the first equation in (\ref{nls}) with $iu$ and $iT_Hu\in L^\infty((0,\infty);E),$ where $T_H$ is defined in the proof of Theorem~\ref{thmmaxdecaylinear}. By (\ref{hypf2'}), $u$ satisfies conservation of charge and we proceed as in the proof of Theorem~\ref{thmmaxdecaylinear}. Hence the claim. See Example~\ref{single2} of Section~\ref{example} for an application.
\end{rmk}

\begin{rmk}
\label{rmkmaxdecaylinear4}
Except for special cases (see for instance Ozawa~\cite{MR91m:35183}), we do not know if the conclusion of Theorem~\ref{thmmaxdecaylinear} still holds when $\|u(t)\|_{L^r(\Omega\cap\{|x-c|<M_0|t|\})}$ is replaced with $\|u(t)\|_{L^r(\Omega\cap\{k'|t|<|x|<k|t|\})}$ in (\ref{thmmaxdecaylinear2}), for some $0<k'<k<\infty.$ The only result known in this direction is the following. Assume that $\Omega=\R^N.$ Then for every $u_0\in L^2(\R^N),$
\begin{gather}
\label{rmkmaxdecaylinear4-1}
\forall k>k'\ge0,\;
\vlim_{t\to\pm\infty}\|e^{it\Delta}u_0\|_{L^2(\{k'|t|<|x|<k|t|\})}=\|\vF^{-1}u_0\|_{L^2(\{k'/2<|x|<k/2\})},
\end{gather}
where $\vF^{-1}$ denotes the inverse Fourier transform in $\R^N.$ Estimate (\ref{rmkmaxdecaylinear4-1}) has been established by Strauss~\cite{297.35062} with $k'=0,$ for any $k>0.$ (see Lemma~p.69 in Strauss~\cite{297.35062} and also Barab~\cite{MR86a:35121}, proof of Lemma~2). From this, we deduce easily that estimate (\ref{rmkmaxdecaylinear4-1}) holds for any $k>k'\ge0.$ It follows that the Strauss's method to establish (\ref{rmkmaxdecaylinear4-1}) uses the Fourier transform and, consequently, does not work when $\Omega\not=\R^N.$ And even for the free operator, we do not know how to adapt his proof when $\Omega\subsetneq\R^N$ is a general domain. Finally, note that the results in Ozawa~\cite{MR91m:35183} are obtained from (\ref{rmkmaxdecaylinear4-1}).
\end{rmk}

\section{Applications}
\label{example}

In this section, we give some examples of nonlinearities for which Theorem~\ref{thmmaxdecaylinear} applies. As is well-known, nonlinear Schrödinger equations enjoy of conservation of a certain energy $E,$ under suitable conditions of the nonlinearity $f.$ In some cases, this implies that $\|\nabla u\|_{L^\infty(\R;L^2(\Omega))}^2\le 2E(u_0)$ $(u$ being a solution of (\ref{nls}) with $u(0)=u_0).$
\begin{exa}
\label{free}
{\bf The free operator} \\
Let $\Omega\subseteq\R^N$ be a nonempty open set and let $(e^{it\Delta})_{t\in\R}$ be the group of isometries generated by $i\Delta$ with the Dirichlet boundary condition on $L^2(\Omega).$ Given $u_0\in H^1_0(\Omega)\setminus\{0\},$ it follows that $u(t)=e^{it\Delta}u_0$ satisfies $iu_t+\Delta u=0$ in $\R\times\Omega.$ Furthermore, the conservation of charge holds and for any $t\in\R,$ $\|\nabla u(t)\|_{L^2(\Omega)}=\|\nabla u_0\|_{L^2(\Omega)}.$ Applying Theorem~\ref{thmmaxdecaylinear} with $f\equiv0,$ we have for any $0<\eps<\|u_0\|_{L^2(\Omega)},$
\begin{gather*}
 \liminf_{t\to\pm\infty}|t|^{N\left(\frac{1}{2}-\frac{1}{r}\right)}\|u(t)\|_{L^r(\Omega\cap\{|x-c|<M_0|t|\})}
 \ge\dfrac{\eps}{|B(0,M_0)|^{\frac{1}{2}-\frac{1}{r}}},
\end{gather*}
for any $c\in\Omega$ and $r\in[2,\infty],$ where
$$
M_0=\frac{2\|u_0\|_{L^2(\Omega)}\|\nabla u_0\|_{L^2(\Omega)}}{\|u_0\|_{L^2(\Omega)}^2-\eps^2}.
\medskip
$$
\end{exa}

\begin{exa}
\label{potential}
{\bf The linear Schrödinger equation with external potential} \\
Let $\Omega\subseteq\R^N$ be an open subset and let $f(u)=-Vu$ where $V\colon\Omega\tends\R$ is a real-valued potential such that $V\in L^p(\Omega)+L^\infty(\Omega),$ for some $p\ge1$ and $p>\frac{N}{2}.$ Then for a given $u_0\in H^1_0(\Omega)\setminus\{0\},$ there exists a unique solution
$$
u\in C_\b(\R;H^1_0(\Omega))\cap C^1_\b(\R;H^{-1}(\Omega))
$$
of (\ref{nls}) such that $u(0)=u_0.$ Moreover, $u$ satisfies conservation of charge and energy $E,$ which is given by
$$
E(u_0)=\dfrac{1}{2}\|\nabla u_0\|_{L^2(\Omega)}^2+\dfrac{1}{2}\dsp\vint_{\Omega}V(x)|u_0(x)|^2dx.
$$
Then it follows from Theorem~\ref{thmmaxdecaylinear} that  for any $0<\eps<\|u_0\|_{L^2(\Omega)},$
\begin{gather*}
 \liminf_{t\to\pm\infty}|t|^{N\left(\frac{1}{2}-\frac{1}{r}\right)}\|u(t)\|_{L^r(\Omega\cap\{|x-c|<M_0|t|\})}
 \ge\dfrac{\eps}{|B(0,M_0)|^{\frac{1}{2}-\frac{1}{r}}},
\end{gather*}
for any $c\in\Omega$ and $r\in[2,\infty],$ where
$$
M_0=\frac{2\|u_0\|_{L^2\Omega)}\|\nabla u\|_{L^\infty(\R;L^2(\Omega))}}{\|u_0\|_{L^2(\Omega)}^2-\eps^2}.
$$
In furthermore $V\ge0,$ then we may choose $M_0$ as
\begin{gather*}
M_0=\frac{2\|u_0\|_{L^2(\Omega)}\sqrt{2E(u_0)}}{\|u_0\|_{L^2(\Omega)}^2-\eps^2}.
\medskip
\end{gather*}
\end{exa}

\begin{exa}
\label{hartree}
{\bf The Hartree type nonlinearity} \\
Let $\Omega=\R^N$ and let $f(u)=-(W*|u|^2)u$ where $W\in L^p(\R^N)+L^\infty(\R^N),$ for some $p\ge1$ and $p>\frac{N}{4},$ and $W^-\in L^q(\R^N)+L^\infty(\R^N),$ for some $q\ge1$ and $q\ge\frac{N}{2}$ $(q>1$ if $N=2).$ We may choose, for example, $W(x)=\mu|x|^{-\gamma}$ with $\mu\in\R\setminus\{0\}$ and $0<\gamma<\min\{N,4\},$ with in addition $\mu>0$ or $\mu<0$ and $0<\gamma<2.$ Then for a given $u_0\in H^1(\R^N)\setminus\{0\},$ there exists a unique solution
$$
u\in C_\b(\R;H^1(\R^N))\cap C^1_\b(\R;H^{-1}(\R^N))
$$
of (\ref{nls}) such that $u(0)=u_0.$ Moreover, $u$ satisfies conservation of charge and energy $E,$ where it is defined by
$$
E(u_0)=\dfrac{1}{2}\|\nabla u_0\|_{L^2(\R^N)}^2+\dfrac{1}{4}\dsp\vint_{\R^N}(W*|u_0|^2)(x)|u_0(x)|^2dx.
$$
Then it follows from Theorem~\ref{thmmaxdecaylinear} that  for any $0<\eps<\|u_0\|_{L^2(\R^N)},$
\begin{gather*}
 \liminf_{t\to\pm\infty}|t|^{N\left(\frac{1}{2}-\frac{1}{r}\right)}\|u(t)\|_{L^r(\{|x-c|<M_0|t|\})}
 \ge\dfrac{\eps}{|B(0,M_0)|^{\frac{1}{2}-\frac{1}{r}}},
\end{gather*}
for any $c\in\R^N$ and $r\in[2,\infty],$ where
$$
M_0=\frac{2\|u_0\|_{L^2(\R^N)}\|\nabla u\|_{L^\infty(\R;L^2(\R^N))}}{\|u_0\|_{L^2(\R^N)}^2-\eps^2}.
$$
In the particular case where $W\ge0$ (or $\mu>0$ if $W(x)=\mu|x|^{-\gamma}),$ we may choose $M_0$ as
\begin{gather*}
M_0=\frac{2\|u_0\|_{L^2(\R^N)}\sqrt{2E(u_0)}}{\|u_0\|_{L^2(\R^N)}^2-\eps^2}.
\end{gather*}
\end{exa}

\begin{exa}
\label{single1}
{\bf Single power interaction in the \boldmath$H^1$\unboldmath--subcritical and attractive case} \\
Let $\Omega\subseteq\R^N$ be an open subset, let $\lambda>0,$ let $0\le\alpha<\frac{4}{N-2}$ $(0\le\alpha<\infty$ if $N\le2)$  and let $f(u)=\lambda|u|^\alpha u.$ The associated energy $E$ is defined by
$$
E(u_0)=\dfrac{1}{2}\|\nabla u_0\|_{L^2(\Omega)}^2
-\dfrac{\lambda}{\alpha+2}\|u_0\|_{L^{\alpha+2}(\Omega)}^{\alpha+2}.
$$
Assume further that one of the following holds. \\
a) $0\le\alpha<\frac{4}{N}.$ \\
b) $\alpha=\frac{4}{N}$ and $u_0\in H^1_0(\Omega)$ with $\|u_0\|_{L^2(\Omega)}$ small enough. \\
c) $\alpha\ge\frac{4}{N}$ and $\|u_0\|_{H^1_0(\Omega)}$ small enough. \\
It follows from Proposition~4.2.3 of Cazenave~\cite{caz1} or Theorem~3.3.5 of Cazenave~\cite{MR2002047}, that for a given $u_0\in H^1_0(\Omega)\setminus\{0\},$ there exists a solution
$$
u\in L^\infty(\R;H^1_0(\Omega))\cap W^{1,\infty}(\R;H^{-1}(\Omega))
$$
of (\ref{nls}) such that $u(0)=u_0$ (see also the proof of Theorem~3.4.3 in Cazenave~\cite{MR2002047}). Moreover, $u$ satisfies conservation of charge and for any $t\in\R,$ $E(u(t))\le E(u_0).$ Note that $u\in C(\R;H^{-1}(\Omega))$ and so $u$ is weakly continuous from $\R$ onto $H^1_0(\Omega).$ Then for any $t\in\R,$ $u(t)\in H^1_0(\Omega)$ and so, $\|u(t)\|_{L^2(\Omega)}$ and $E(u(t))$ are well defined for any $t\in\R.$ Furthermore, it follows from Remark~\ref{rmkmaxdecaylinear3} that $u\in C(\R;L^2(\Omega)).$ By Theorem~\ref{thmmaxdecaylinear}, we have for any $0<\eps<\|u_0\|_{L^2(\Omega)},$
\begin{gather*}
 \liminf_{t\to\pm\infty}|t|^{N\left(\frac{1}{2}-\frac{1}{r}\right)}\|u(t)\|_{L^r(\Omega\cap\{|x-c|<M_0|t|\})}
 \ge\dfrac{\eps}{|B(0,M_0)|^{\frac{1}{2}-\frac{1}{r}}},
\end{gather*}
for any $c\in\Omega$ and $r\in[2,\infty],$ where
$$
M_0=\frac{2\|u_0\|_{L^2(\Omega)}\|\nabla u\|_{L^\infty(\R;L^2(\Omega))}}{\|u_0\|_{L^2(\Omega)}^2-\eps^2}.
$$
For the repulsive case $\lambda<0,$ see Example~\ref{single2}.
\medskip
\end{exa}

For more details about global existence and boundness in the energy space in Examples~\ref{hartree}--\ref{single1}, see for example Bourgain~\cite{MR2000h:35147}, Cazenave~\cite{caz1,MR2002047}, Ginibre~\cite{ginibre}, C.~Sulem and P.-L.~Sulem~\cite{MR2000f:35139} and the references therein.

\begin{exa}
\label{single2}
{\bf Large nonlinearity -- The repulsive case} \\
Let $\Omega\subseteq\R^N$ be an open subset and let $f(u)=-\lambda|u|^\alpha u$ with $\lambda>0$ and $0\le\alpha<\infty.$ We set $E=H^1_0(\Omega)\cap L^{\alpha+2}(\Omega)$ and $E^*=H^{-1}(\Omega)+L^\frac{\alpha+2}{\alpha+1}(\Omega).$ Since $\vD(\Omega)$ is dense in both $H^1_0(\Omega)$ and $L^{\alpha+2}(\Omega),$ then we have from Lemma~2.3.1 and Theorem~2.7.1 in Bergh and Löfström~\cite{bl} that  $E$ and $E^*$ are reflexive separable Banach spaces and that $E^*$ is the topological dual of $E.$ Given $u_0\in H^1_0(\Omega)\setminus\{0\},$ there exists a solution $u\in L^\infty(\R;E)\cap W^{1,\infty}(\R;E^*)$ of (\ref{nls}) with initial value $u_0.$ Moreover, $u$ satisfies conservation of charge and for any $t\in\R,$ $E(u(t))\le E(u_0),$ where 
\begin{gather}
\label{energy}
E(u_0)=\dfrac{1}{2}\|\nabla u_0\|_{L^2(\Omega)}^2
+\dfrac{\lambda}{\alpha+2}\|u_0\|_{L^{\alpha+2}(\Omega)}^{\alpha+2}.
\end{gather}
See Strauss~\cite{MR80b:35090}; see also Cazenave~\cite{caz1,MR2002047}, Section~9.4. Note that $u\in C(\R;E^*)$ and so $u$ is weakly continuous from $\R$ onto $H^1_0(\Omega)$ and from $\R$ onto $L^{\alpha+2}(\Omega).$ Then for any $t\in\R,$ $u(t)\in E$ and so, $\|u(t)\|_{L^2(\Omega)}$ and $E(u(t))$ are well defined for any $t\in\R.$ Furthermore, we have by Remark~\ref{rmkmaxdecaylinear3} that $u\in C(\R;L^2(\Omega)).$ It is clear that $f$ satisfies (\ref{hypf1'})--(\ref{hypf3'}) and so it follows from Theorem~\ref{thmmaxdecaylinear} and Remark~\ref{rmkmaxdecaylinear3} that for any $0<\eps<\|u_0\|_{L^2(\Omega)},$
\begin{gather*}
 \liminf_{t\to\pm\infty}|t|^{N\left(\frac{1}{2}-\frac{1}{r}\right)}\|u(t)\|_{L^r(\Omega\cap\{|x-c|<M_0|t|\})}
 \ge\dfrac{\eps}{|B(0,M_0)|^{\frac{1}{2}-\frac{1}{r}}},
\end{gather*}
for any $c\in\Omega$ and $r\in[2,\infty],$ where
$$
M_0=\frac{2\|u_0\|_{L^2(\Omega)}\sqrt{2E(u_0)}}{\|u_0\|_{L^2(\Omega)}^2-\eps^2}.
$$
Note that when $N\le2,$ or $N\ge3$ and $\alpha\le\frac{4}{N-2},$ we have that $E=H^1_0(\Omega)$ and $E^*=H^{-1}(\Omega).$
\end{exa}

\begin{rmk}
Notice that Examples~\ref{potential}, \ref{hartree} and \ref{single2} cover the results of Hayashi and Ozawa in \cite{MR91d:35035}.
\end{rmk}

\begin{rmk}
Here are some comments about uniqueness and smoothness in the case of the single power of interaction.
\begin{enumerate}
\item
In Examples~\ref{single1} and \ref{single2}, if $N=1$ then the solution $u$ is unique and $u\in C_\b(\R;H^1_0(\Omega))\cap C^1_\b(\R;H^{-1}(\Omega)).$ Furthermore, conservation of energy holds. See for example Theorem~4.4.1 in Cazenave~\cite{caz1} or Theorem~3.5.1 in Cazenave~\cite{MR2002047}.
\item
In Examples~\ref{single1} and \ref{single2}, if $N=2$ and if furthermore $\Omega=\R^N$ or $\alpha\le\frac{4}{N},$ then the solution $u$ is unique and $u\in C_\b(\R;H^1_0(\Omega))\cap C^1_\b(\R;H^{-1}(\Omega)).$ Furthermore, conservation of energy holds. See for example Theorem~4.5.1 in Cazenave~\cite{caz1} or Theorem~3.6.1 in Cazenave~\cite{MR2002047}.
\item
In Examples~\ref{single1} and \ref{single2}, if $N\ge3$ and if furthermore $\Omega=\R^N$ and $\alpha<\frac{4}{N-2},$ then the solution $u$ is unique and $u\in C_\b(\R;H^1(\R^N))\cap C^1_\b(\R;H^{-1}(\R^N)).$ Furthermore, conservation of energy holds. See for example Theorem~4.3.1 in Cazenave~\cite{caz1,MR2002047} or Corollary~4.3.3 in Cazenave~\cite{MR2002047}.
\item
For exterior domains, we have the following result. Let $N\ge2$ be an integer,  $C\subset\R^N$ be a nonempty star-shaped open bounded subset with smooth boundary, $\Omega=\R^N\setminus\ovl C$ and let $f(u)=\lambda(1+|u|^2)^\frac{\alpha}{2}u$ with $\lambda>0$ and $0\le\alpha<\frac{4}{N}.$ If $N\ge5$ then assume further that $\alpha<\frac{2}{N-2}.$ Using the result of Burq, Gérard and Tzvetkov~\cite{MR2068304}, we obtain that for any $u_0\in H^1_0(\Omega)\setminus\{0\},$ there exists a unique solution
$
u\in C_\b(\R;H^1_0(\Omega))\cap C^1_\b(\R;H^{-1}(\Omega))
$
of (\ref{nls}) such that $u(0)=u_0,$ satisfying conservation of charge and energy $E,$ where
$$
E(u_0)=\dfrac{1}{2}\|\nabla u_0\|_{L^2(\Omega)}^2
-\dfrac{\lambda}{\alpha+2}\vint_\Omega\left((1+|u_0(x)|^2)^\frac{\alpha+2}{2}-1\right)dx.
$$
Finally, Theorem~\ref{thmmaxdecaylinear} applies.
\end{enumerate}
\end{rmk}

\begin{rmk}
We may mix some nonlinearities by taking
$$
f(u)=\lambda|u|^\alpha u+\mu|u|^\gamma u+Vu+(W*|u|^2)u, \mbox{ in } \Omega=\R^N,
$$
or
$$
f(u)=\lambda|u|^\alpha u+\mu|u|^\gamma u+Vu, \mbox{ in } \Omega\subseteq\R^N,
$$
for some constants $\lambda,\mu,\alpha,\gamma\in\R$ and some real-valued functions $V$ and $W.$ See for example Cazenave~\cite{caz1} (Sections 4.3--4.5 and 6.5) or Cazenave~\cite{MR2002047} (Sections~3, 4.1--4.5 and 6.8) and the references therein.
\end{rmk}

\noindent
{\large\bf Acknowledgments} \\
The author is grateful to Professor Thierry Cazenave for stimulating remarks and for pointing out to him several imprecisions in the first version of the paper.

\baselineskip .0cm

\bibliographystyle{abbrv}
\bibliography{BiblioPaper4}

\end{document}